# A Guided Tour in the Topos of Graphs

Sebastiano Vigna*

**Abstract**

In this paper we survey the fundamental constructions of a presheaf topos in the case of the elementary topos of graphs. We prove that the transition graphs of nondeterministic automata (a.k.a. labelled transition systems) are the separated presheaves for the double negation topology, and obtain as an application that their category is a quasitopos.

## 1 Introduction

Sheaf topoi are usually associated to continuous mathematics, such as differential or algebraic geometry. Nonetheless, several (pre)sheaf topoi with simple (even finite!) base categories describe fundamental mathematical objects such as trees and graphs. It has been suggested by Lawvere [5] that these simple topoi have a rich combinatorial structure. Nonetheless, very few results are known about them, and their internal workings seem to be mostly unexplored.

This note takes a brief tour through the topos of graphs, a presheaf topos having as base the category with two parallel arrows [5, 12, 6] (here and in the following we assume that our graphs are directed, and that they may possess parallel arcs and self-loops; note that elsewhere these graphs have been called *irreflexive*, as they do not possess an identity loop at each vertex). We introduce the basic definitions of topos theory and show their meaning in the topos of graphs, always trying to describe intuitively the essence of each construction. We assume some knowledge of category theory, in particular the notion of cartesian closedness. Good references are [7], [2] and [8]. Moreover, an elementary introduction to topoi can be found in [6].

Our purpose is twofold: on one side, our survey on the topos of graphs provides a detailed, complete and self-contained example of "combinatorial topos", which can be extremely useful to anyone wanting to get acquainted with topos-theoretical definitions; on the other side, we show that several constructions of the topos correspond to known constructions in graph theory (such as the induced subgraph construction), and by analyzing the Lawvere topologies of the topos we obtain a characterization of the category of automata as a category of separated objects of a topos; this allows us to prove immediately that automata (i.e., labelled transition systems) and other mathematical objects, such as simple finite undirected graphs and tolerance spaces, form a quasitopos.

---

*Dipartimento di Scienze dell'Informazione, Università di Milano, Via Comelico 39/41, I-20135 Milano MI, Italy, Fax: +39-2-55006276. e_mail: `vigna@dsi.unimi.it`

Our treatment is less elementary than the one given in [6] (yet we do have some minor overlaps, in particular in the detailed description of Ω given in Section 3), as it assumes a certain knowledge of category theory; on the other hand, we strive to minimize the number of notions introduced, so to enlarge the audience. This sometimes keeps us from making some interesting remarks that, however, will be evident to the specialists, by whom we hope to be forgiven.

## 2 Basic definitions

We begin with the fundamental definition of elementary topos:

**Definition 1** An *elementary topos* is a cartesian closed category with a *subobject classifier*, i.e., an object $\Omega$ with an arrow $\top : 1 \to \Omega$ such that for every mono $m : S \hookrightarrow X$ there is a unique arrow $\chi_m : X \to \Omega$ which makes the following diagram a pullback:

$$\begin{array}{ccc} S & \longrightarrow & 1 \\ m \downarrow & & \downarrow \top \\ X & \xrightarrow{\chi_m} & \Omega \end{array}$$

the arrow $\chi_m$ is the *classifying*, or *characteristic* arrow for $m$.

$\chi_m$ generalizes the well-known characteristic map from set theory (in the topos of sets, $\Omega = \{0, 1\}$ and the pullback exhibits the well-known correspondence between subsets and characteristic functions). The subobject classifier can also be seen as an object of *truth values*: the arrow $\top$ is called the *true* map, while the *false* map $\bot : 1 \to \Omega$ is the classifier of the unique map $0 \to 1$ (with 0 and 1 we denote the initial and terminal objects, respectively; it can be shown that every topos has all finite colimits).

Here we shall deal with a particular class of elementary topoi, i.e., those constructed by means of *presheaf categories*: given a small category $\mathcal{C}$, the category $\mathbf{Sets}^{\mathcal{C}^{op}}$ of functors $\mathcal{C}^{op} \to \mathbf{Sets}$ and natural transformations turns out to be a topos. In fact, this statement is true for any generalized functor category $\mathcal{E}^{\mathcal{C}^{op}}$, as long as $\mathcal{E}$ is a topos and $\mathcal{C}$ is a category in $\mathcal{E}$ (for the precise definition, see [8]); this technique gives a way of building new topoi from given ones. For instance, presheaves of finite sets with a finite base category form a topos.

Let now consider the category $\Gamma$ represented in the following picture:

$$N \underset{t}{\overset{s}{\rightrightarrows}} A$$

The presheaf category $\mathcal{G} = \mathbf{Sets}^{\Gamma^{op}}$ is the topos of graphs. Each presheaf $G$ in $\mathcal{G}$ is given by a set $G(N)$, the set of *nodes*, and a set $G(A)$ of *arcs*. The arrows $s, t$ are mapped to functions $G(s), G(t) : G(A) \to G(N)$ which assign to each arc its *source* and *target*. The reader can easily check that natural transformations are the classical graph morphisms (i.e., maps of nodes and arcs which preserve the source/target assignments). Analogously, $\mathbf{FinSets}^{\Gamma^{op}}$ is the topos of finite graphs (our results are true both in the finite and in the general case).

# 3  Representables, limits, exponentials and the classifier

Presheaf categories possess standard constructions for (co)limits, for $\Omega$ and for the exponentials. We are going to discuss these constructions, giving the general definitions and applying them to $\mathcal{G}$. We do not give proofs, for which the reader can refer to the abovementioned texts.

The base of our constructions is given by the representable functors. There is one representable $\mathcal{C}(-, X)$ for each object $X$ of $\mathcal{C}$, which associates to $X$ the set of morphisms into $X$. The definition on arrows is given by composition.

In $\mathcal{G}$, the representable functors are $N = \Gamma(-, N)$ and $A = \Gamma(-, A)$. Since $\Gamma(N, N) = \{N\}$ and $\Gamma(A, N) = \varnothing$, the representable $N$ has just one node, called $N$ (note that, as usual, when no confusion is possible we write $X$ for $\mathbf{1}_X$). The representable $A$ is the one-arrow graph pictured below:

$$s \xrightarrow{A} t \;,$$

because $\Gamma(N, A) = \{s, t\}$ and $\Gamma(A, A) = \{A\}$.

One can think of representable functors as a representation in the "real world" (the topos) of "concepts" (objects of the base category): in our case, the representable associated to the object $N$, whose image is the set of nodes, is a single node, while the representable associated to the object $A$, whose image is the set of arcs, is a single arc.

It can be shown that the representables provide a set of generators for the topos (i.e., for every $f, g : X \to Y$ we have $f = g$ iff $f \circ h = g \circ h$ for each representable $R$ and each arrow $h : R \to X$). Thus, two arrows are different iff they are different over some representable. Moreover, every object $X$ of a topos is built, in a (precise) sense, by glueing a suitable set of representables, because $X$ is the colimit of the diagram $r : R \to X$, where $R$ ranges over the representables and $r$ ranges over the morphisms $R \to X$.

In a presheaf topos, sums and products are computed locally, i.e., the (co)limit $L$ of any functor $F : \mathcal{D} \to \mathbf{Sets}^{\mathcal{C}^{op}}$ is given by

$$L(X) = (\text{co}) \lim F(-)(X),$$

where $F(-)(X) : \mathcal{D} \to \mathbf{Sets}$ is the functor obtained by fixing the second coordinate of $F$ (which can be seen, by cartesian closedness, as a functor $\mathcal{D} \times \mathcal{C}^{op} \to \mathbf{Sets}$)[1]. Using this fact, it is easy to check that 0 is the empty graph (i.e., $0(N) = 0(A) = \varnothing$), while 1 is a self-loop (because $1(N) = 1(A) = \{*\}$). Moreover, the product $G \times H$ of two graphs has node set $G(N) \times H(N)$ and arc set $G(A) \times H(A)$; in other words, we put an arc between $\langle x, y \rangle$ and $\langle x', y' \rangle$ for each pair of arcs between from $x$ to $x'$ (in $G$) and from $y$ to $y'$ (in $H$).

In general, the exponential $Y^X$ in the category $\mathbf{Sets}^{\mathcal{C}^{op}}$ is defined by

$$Y^X(Z) = \mathbf{Nat}(\mathcal{C}(-, Z) \times X, Y) = \mathbf{Sets}^{\mathcal{C}^{op}}(\mathcal{C}(-, Z) \times X, Y),$$

and the definition on arrows follows by naturality ($\mathbf{Nat}(-, -)$ is, of course, the set of natural transformations between two functors). In our case, we have

$$H^G(N) = \mathcal{G}(N \times G, H).$$

---

[1] More generally, the forgetful functor $U : \mathbf{Sets}^{\mathcal{C}^{op}} \to \mathbf{Sets}/|\mathcal{C}|$ creates (co)limits.

Note that the effect of multiplying $G$ by $N$ is just to strip away all of its arcs. Thus, the nodes of $H^G$ are the functions $G(N) \to H(N)$. For what matters arcs, we have

$$H^G(A) = \mathcal{G}(A \times G, H).$$

The graph $A \times G$ is built by taking two copies of $G(N)$, and attaching each arc in $G(A)$ in such a way that its source is in the first copy, while its target is in the second copy. Equivalently, we can think of taking $G$, adding a copy of $G(N)$, and then moving the target of each arc to the second copy.

Each morphism $A \times G \to H$ is an arc of $H^G$, having as source (target) the map induced by the first (second) copy of $G(N)$ in $A \times G$. Intuitively, an arc of $H^G$ is a "relaxed" morphism $G \to H$ which maps arcs with a common source (target) to arcs with a common source (target). Of course, such a morphism does not correspond, in general, to a morphism $G \to H$: in fact, it is as farther from being such a morphism as its source and target are "different": when they are equal, the arc is a self-loop, and it represents a morphism $G \to H$ (as it should, since the *global elements* $1 \to H^G$ are in bijection with the morphisms $G \to H$).

In the next picture, we show a worked out example. The second copy of $G(N)$ is denoted by overlines, and we described the maps $\{a, b, c\} \to \{0, 1\}$ (or $\{\bar{a}, \bar{b}, \bar{c}\} \to \{0, 1\}$) by listing in order the number assigned to each element.

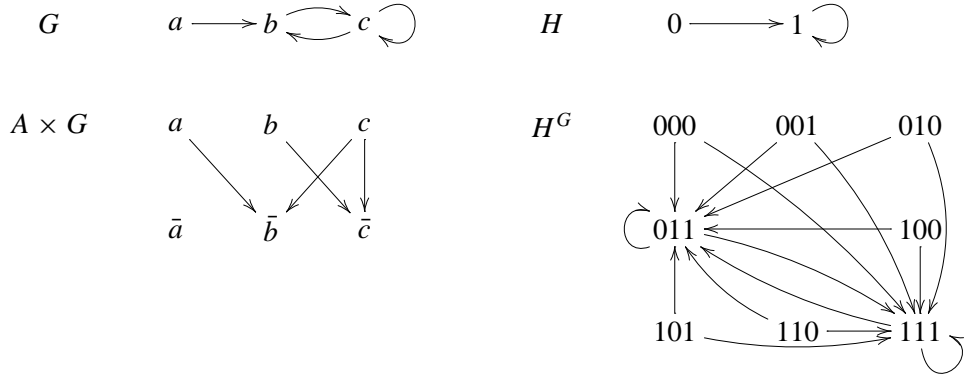

Finally, the classifier $\Omega$ can be built by considering subobjects of the representables. In general, for $\mathbf{Sets}^{\mathcal{C}^{op}}$ we have

$$\Omega(X) = \{S \mid S \text{ is a subobject of } \mathcal{C}(-, X)\},$$

and the morphism $\top : 1 \to \Omega$ is given by $\top(X) = X$ (i.e., its value on $X$ is the representable on $X$).

In our case, the representable $N$ has just two subobjects, $0_N$ and $N$, which provide the node set of $\Omega$, while the representable $A$ has five subobjects, namely $0_A$, $A$, $s$, $t$ and $\binom{s}{t}$, which provide the arc set[2]. This gives us the following graph:

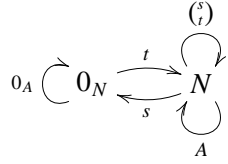

where the "true" arrow $\top : 1 \to \Omega$ sends the only node and arc of 1 to $N$ and $A$, respectively.

Given a subgraph $m : S \hookrightarrow G$, the classifying map $\chi_m : G \to \Omega$ works as follows:

- nodes which are not in $S$ are mapped to $0_N$;
- nodes which are in $S$ are mapped to $N$;
- if an arc is not in $S$, we have four possibilities:
    - arcs whose source and target are not in $S$ are mapped to $0_A$;
    - arcs whose source is in $S$, but the target is not, are mapped to $s$;
    - arcs whose target is in $S$, but the source is not, are mapped to $t$;
    - arcs having both source and target in $S$ are mapped to $\binom{s}{t}$;
- finally, an arc in $S$ is mapped to $A$.

In other words, the possibilities for nodes are just two (a node may, or may not, be in a certain subgraph); however, the situation is much more varied for an arc, to which five different arcs of $\Omega$ can be assigned, depending also on the assignment of its source and target.

As we already remarked, $\Omega$ is the object of truth values of the topos. Logical operators on these truth values can be defined as follows: the conjunction $\wedge : \Omega \times \Omega \to \Omega$ is the characteristic map of the subobject $\langle \top, \top \rangle : 1 \to \Omega \times \Omega$; the negation $\neg : \Omega \to \Omega$ is the characteristic map of $\bot : 1 \to \Omega$.

It is easy to check that on the nodes $\wedge$ and $\neg$ work exactly like in **Sets** (as we already mentioned, nodes have only two truth values). On the arcs the situation is different: $\wedge$ "reduces" the truth value of every pair of arcs of $\Omega$ to the minimum (the order being $0_A < s, t < \binom{s}{t} < A$); for instance, $s \wedge t = 0_A$, while $\binom{s}{t} \wedge s = s$. The negation map, on the other hand, satisfies $\neg 0_A = A$ and $\neg A = \neg \binom{s}{t} = 0_A$ (so that it "forgets" whether an arc was in a subgraph or not, as long as its source and target were).

## 4 Topologies on $\Omega$

Another fundamental source of new topoi is given by *topologies*: a topology allows to set out certain objects, called *sheaves*, which form a new topos (in fact, sheaves were the motivating example for topos theory).

---

[2]We are using a number of common abbreviations: for every graph $G$, $0_G$ is the unique morphism $0 \to G$. The maps $s, t : N \to A$ are defined in the obvious way, i.e., $s(N) = s$. The map $\binom{s}{t} : N + N \to A$ is induced by the universal property of the sum. Often we do not distinguish explicitly the node and arc component of a graph morphism.

There are many ways to define a topology on a topos. Here we will follow Lawvere's idea of an elementary description, given by a map $j : \Omega \to \Omega$ which satisfies a simple set of axioms.

**Definition 2** A *topology* on a topos is a morphism $j : \Omega \to \Omega$ such that

$$j \circ \top = \top \tag{1}$$
$$j \circ j = j \tag{2}$$
$$j \circ \wedge = \wedge \circ (j \times j). \tag{3}$$

In other words, $j$ preserves conjunction and truth, and it is idempotent. The original idea behind a Lawvere topology is an (internal) generalization of a Grothendieck topology (which in turn is a categorical generalization of a classical topology); nonetheless, the definition also applies in situations where continuity is apparently absent. In this case, $j$ becomes a combinatorial, rather than a topological object.

In order to see this fact more clearly, we classify the nontrivial topologies on $\mathcal{G}$. We remark that the map $\neg\neg : \Omega \to \Omega$ is always a topology, called the *double negation* topology, and that two trivial topologies are always available: the identity on $\Omega$ and the morphism $\Omega \xrightarrow{!_\Omega} 1 \xrightarrow{\top} \Omega$.

**Theorem 1** There is exactly one nontrivial topology on $\mathcal{G}$ besides the double negation topology.

**Proof.** The proof is elementary (in fact, a case-by-case analysis). First of all note that $j(N) = N$ and $j(A) = A$ necessarily by (1).

If $j(0_N) = 0_N$, the only choice we have is to send $\binom{s}{t}$ to itself (which implies $j = \mathbf{1}_\Omega$) or to $A$ (the image of $0_A$, $t$ and $s$ is forced by our choices for the nodes); in the last case, we obtain the double negation topology, because $\neg\neg : \Omega \to \Omega$ fixes all arcs (and nodes), except for $\binom{s}{t}$, which is sent to $A$.

If $j(0_N) = N$, then $j(0_A) = A$ implies $j = \top \circ !_\Omega$; indeed, by (3) we have

$$A = j(0_A) = j\left(\binom{s}{t} \wedge 0_A\right) = j\left(\binom{s}{t}\right) \wedge j(0_A) = j\left(\binom{s}{t}\right) \wedge A = j\left(\binom{s}{t}\right),$$

and since by (2) $j(t)$ and $j(s)$ must be arcs fixed by $j$, we have $j(s) = j(t) = A$.

Thus, we must assume $j(0_A) = \binom{s}{t}$ (which implies by (3) $j(s) = j(t) = \binom{s}{t}$). There are no other remaining cases. ∎

Another possible statement of the last theorem is that the lattice of topologies on $\Omega$ is the four-elements boolean algebra. The interested reader can check that the only nontrivial topology besides double negation is the so-called *closed topology* on the global element $\binom{s}{t} : 1 \to \Omega$, which is defined as $- \vee \binom{s}{t}$ [8]. Since this implies that the double negation topology is open, we shall call the only nontrivial topology different from $\neg\neg$ *the* closed topology.

## 5 Closure and density

**Definition 3** Given a subobject $m : S \hookrightarrow G$, with characteristic map $\chi_m : G \to \Omega$, the closure of $S$ in $G$ (with respect to the topology $j$) is the subobject $\bar{S}$ classified by $j \circ \chi_m$; $S$ is said to be *dense* if its closure is equal to $G$.

Intuitively, composition with $j$ "increases" the truth level of $\chi_m$, thus describing a bigger object. The axioms on $j$ guarantee that the behaviour of this closure operation is reasonable, i.e., $S \hookrightarrow \bar{S}$, $\bar{\bar{S}} = \bar{S}$.

The next two theorems show that the topologies on $\Omega$ generate closure operations and dense subobjects which correspond to well-known constructions from graph theory. Recall that a *spanning subgraph* (or *partial graph*) of $G$ is a subgraph containing all the nodes of $G$, while an *induced subgraph* is a subgraph $H$ containing all the arcs of $G$ connecting nodes of $H$.

**Theorem 2** The closure associated to the closed topology adds to a subgraph $S$ all the nodes of the graph, i.e., $\bar{S}$ is the spanning subgraph generated by $S$. The dense subobjects of $G$ are those subgraphs which include all the *arcs* of $G$. In particular, there is a *minimum* dense subobject, which is the arc set of $G$.

**Proof.** Let $m : S \hookrightarrow G$ be a subgraph of $G$, and $\chi_m$ its characteristic map. The composition with the closed topology has the effect of adding all the nodes of $G$ to $S$, because nodes previously mapped to $0_N$ will be mapped to $N$; moreover, no arc will be added to $S$, because no arc of $\Omega$ is sent to $A$ by the closed topology, except for $A$ itself. This means that a subgraph $S$ needs exactly to include all the arcs of $G$ in order to be dense. ∎

**Theorem 3** The closure associated to the double negation topology adds to a subgraph $H$ all the arrows which have source and target in $S$, i.e., $\bar{S}$ is the induced subgraph generated by $S$. The dense subobjects of $G$ are those subgraphs which include all the *nodes* of $G$. In particular, there is a *minimum* dense subobject, which is the node set of $G$.

**Proof.** Let $m : S \hookrightarrow G$ be a subgraph of $G$, and $\chi_m$ its characteristic map. The composition with $\neg\neg$ has the effect of adding to $S$ all arcs of $G$ lying between nodes of $S$; in fact, such arcs are sent by $\chi_m$ to $\binom{s}{t}$, and this arc is in turn sent to $A$ by $\neg\neg$. This means that a subgraph $S$ needs exactly to include all the nodes of $G$ in order to be dense. ∎

## 6 Sheaves and separated objects

**Definition 4** An object $X$ of a topos with topology $j$ is said to be *j-separated* if for every object $Y$, every $j$-dense subobject $m : S \hookrightarrow Y$, and every morphism $f : S \to X$ there is at most one factorization $g : Y \to X$ of $f$ through $m$:

$$\begin{array}{ccc} S & & \\ {\scriptstyle m}\downarrow & \searrow^{f} & \\ Y & \xrightarrow{g} & X \end{array}$$

An object is *j-complete* such a factorization always exists. A *j-sheaf* is an object which is $j$-separated and $j$-complete[3].

---
[3] We shall usually understand the reference to $j$ if it is clear from the context.

In other words, *a morphism into a separated object X is completely determined by its restriction to a dense subobject of the domain*. If the object $X$ is also complete (i.e., a sheaf), every morphism defined on a dense subobject $S \hookrightarrow Y$ into $X$ extends to a unique morphism $Y \to X$. Sheaves are important because their full subcategory is again a topos.

Getting back to $\mathcal{G}$, it is easy to check that separated graphs and sheaves for the closed topology are rather trivial: sheaves are graphs with just one node, and their topos is equivalent to **Sets** (via the functor $(-)(A) : \mathcal{G} \to$ **Sets**); the only separated, noncomplete object is the empty graph. But in the case of the double negation topology, something radically different happens:

**Theorem 4** The objects of $\mathcal{G}$ which are separated for the double negation topology are exactly the graphs without parallel arcs. The sheaves are exactly the complete graphs (with self-loops), and they form a topos equivalent to **Sets**.

**Proof.** Consider the monomorphism $\binom{s}{t} : N + N \to A$, and a separated graph $G$ with two arcs $a$ and $b$ between the nodes $x$ and $y$. The morphisms $a, b : A \to G$ defined in the obvious way are both extensions of $\binom{x}{y} : N + N \to G$; by separateness, this implies $a = b$. On the other hand, if $G$ is complete for any pair of nodes $x$ and $y$ an extension $a : A \to G$ of $\binom{x}{y} : N + N \to G$ determines an arc between $x$ and $y$. Thus separated (complete) graphs have at most (least) one arc between each pair of nodes.

On the other hand, if $m : S \hookrightarrow H$ is dense, by Theorem 3 the inclusion $m$ has to be a bijection on the nodes; thus, every extension of a map $f : S \to G$ to $H$ has to be defined on nodes as $f_N \circ (s_N)^{-1}$. If $G$ has no parallel arcs, then the extension (if it exists) is also uniquely determined on the arcs, and if $G$ has an arc between each pair of nodes an extension always exists. The last statement is immediate, the equivalence being induced by the restriction of the functor $(-)(N) : \mathcal{G} \to$ **Sets**. ∎

The last theorem suggests that a reasonable name for graphs without parallel arcs should be "separated" (the separated objects for the closed topology being uninteresting). Note that, by one of the mysterious coincidences of mathematics, classical complete graphs are complete in the topos-theoretical sense.

## 7 Arc labellings and transition systems

We now want to apply our separateness results. Consider an alphabet $\Sigma$ and the graph $\Sigma$ (we are slightly abusing the notation) defined by $\Sigma(N) = \{*\}$, $\Sigma(A) = \Sigma$ (the arrows are forced by the universal property of $\{*\}$). The slice category $\mathcal{G}/\Sigma$ is formed in the usual way, by taking as objects the morphisms $g : G \to \Sigma$, and as arrows the commuting triangles

$$\begin{array}{ccc} G & \xrightarrow{a} & H \\ {}_g\searrow & & \swarrow_h \\ & \Sigma & \end{array}$$

The fundamental theorem of topos theory says that every slice of a topos is again a topos, and this fact of course applies to $\mathcal{G}/\Sigma$. Indeed, slicing is a most important source of new topoi.

Since a morphism $g : G \to \Sigma$ must send all nodes of $G$ to $*$, and can choose freely an element of $\Sigma$ for each arc of $G$, $\mathcal{G}/\Sigma$ contains the graphs arc-labelled on $\Sigma$. The subobject classifier is essentially the same, but we must add a copy of each arrow for each element of $\Sigma$ (because morphisms in $\mathcal{G}/\Sigma$ have to preserve labels). Technically, the new subobject classifier is given by the first projection $\Sigma \times \Omega \to \Sigma$ [8].

Our interest is in the transition graphs of (nondeterministic) automata[4], also known as *labelled transition systems*, whose category has been widely studied in computer science, as it is used for the semantics of process algebras (see, for instance, [9] and [4]).

Not all objects of $\mathcal{G}/\Sigma$ are transition systems; indeed, transition systems are exactly those graphs satisfying the following condition: for each pair of nodes (states) $x$, $y$ there is at most one arc (transition) labelled by $\alpha$ with $x$ as source and $y$ as target. Thus, the following graph

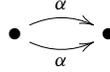

is not a transition system.

It is easy to check that the proof of Theorem 4 works also for the slice topos $\mathcal{G}/\Sigma$; only, this time we must require that there are no parallel arcs with the same label. Indeed, the double negation topology is identical, except that it works labelwise. This yields the following

**Theorem 5** Transition systems labelled on $\Sigma$ are exactly the objects of $\mathcal{G}/\Sigma$ which are separated for the double negation topology.

Moreover, we can now easily prove that

**Theorem 6** Transition systems labelled on $\Sigma$ form a quasitopos.

The result is immediate by Theorem 5 and by [3], where it is shown that separated objects form a quasitopos. Since the category of automata and the category of their transition graphs (with graph morphisms) are equivalent, we have also the following

**Corollary 1** Automata on $\Sigma$ form a quasitopos.

A deep analysis of quasitopoi can be found in [12], where several of their properties, such as local cartesian closedness, finite completeness, representability of relations and so on, are proved. We just note that the essential difference between a topos and a quasitopos is that that $\Omega$ does not classify *all* subobjects, but rather those subobjects defined by *strong monomorphisms*. A monomorphism $m : X \to Y$ is said to be *strong* iff every commutative square

$$\begin{array}{ccc} X' & \xrightarrow{e} & Y' \\ f \downarrow & \overset{d}{\swarrow} & \downarrow g \\ X & \xrightarrow{m} & Y \end{array}$$

with $e$ an epimorphism has a *diagonal $d$* as indicated (this means that the triangles commutes, so $d$ is unique). In the case of transition systems, we have our last

---
[4]By an *automaton* on an alphabet $\Sigma$ we mean a function $\delta : \Sigma \times X \to 2^X$; a morphism between automata $\delta : \Sigma \times X \to 2^X$ and $\delta' : \Sigma \times X' \to 2^{X'}$ is a function $f : X \to X'$ such that $f(\delta(\alpha, x)) \subseteq \delta'(\alpha, f(x))$.

**Theorem 7** The strong monomorphism are those defining induced subgraphs, i.e., a strong subobject of a transition system is defined by the removal of a subset of states.

**Proof.** The necessity of this condition can be seen in the following diagram, where $X$ has an arc $a$ between $x$ and $y$ and $m : S \hookrightarrow X$ contains $x$ and $y$, but not $a$:

$$\begin{array}{ccc} N + N & \xrightarrow{\binom{s}{t}} & A \\ \binom{x}{y} \downarrow & & \downarrow a \\ S & \xhookrightarrow{m} & X \end{array}$$

It is clear that no diagonal exists. Sufficiency can be proved by using the fact that (the monomorphism of) an induced subgraph is always an isomorphism on its image. ∎

## 8  Conclusions

We hope that the reader appreciated the rich internal structure of (the topologies on) $\mathcal{G}$. Of course, other topoi of graphs are also of interest: in particular, by adding an arrow $\varepsilon : A \to N$ to $\Gamma$, satisfying the simple equations $s \circ \varepsilon = t \circ \varepsilon = \mathbf{1}$, we obtain the topos of *reflexive* graphs [5], in which every node has an assigned self-loop which is preserved by morphisms (in fact, nodes should be more correctly *identified* with such self-loops). This does not alter in an essential way the combinatorial structure of a single graph, but now morphisms can collapse arcs (the almost standard name in the graph-theoretical literature is *degenerate map*), so all constructions inside the topos are influenced—one can easily check that the product is now another classical operation on graphs. Shrimpton [11] pursues a most interesting study of the structure of the "group-graph" of automorphisms in this new topos.

Both (reflexive and irreflexive) topoi can be made "undirected" by augmenting $\Gamma$ with a *symmetry* arrow, i.e., an involution on $A$ satisfying obvious equations w.r.t. $s$ and $t$ (and possibly $\varepsilon$). The resulting graphs are very similar to classic undirected graphs, but with a distinguished feature: there are self-loops that are fixed by the symmetry and self-loops that are not, something which is not expressible in terms of undirected graphs. This fact has a strong impact on the combinatorial structure of coverings, as discussed in [1].

However, in all cases graphs that are separated for the double negation topology are graphs without parallel arcs, which reinforces the idea that "separated" should be the right word for such graphs (no matter whether they are also undirected or reflexive). In other words, Theorem 4 applies also to the topoi discussed in this section, so several categories of standard mathematical structures (with their standard morphisms) turn out automatically to be quasitopoi: just to name a few, simple finite undirected (*schlicht*) graphs (but with self-loops allowed), tolerance spaces [10] and (of course) binary endorelations, a.k.a. digraphs without parallel arcs.